\documentclass[amstex,12pt]{article}

\usepackage{graphicx}
\usepackage{amsfonts}
\usepackage{amssymb}
\usepackage{amsmath}
\usepackage{latexsym}
\usepackage{enumerate}

\newtheorem{thm}{Theorem}[section]
\newtheorem{lm}[thm]{Lemma}

\newtheorem{co}[thm]{Corollary}

\newtheorem{ex}[thm]{Example}
\newtheorem{df}[thm]{Definition}

\begin{document}

\author{Jin-ichi Itoh and Costin V\^\i lcu}
\title{Orientable cut locus structures on graphs}
\maketitle

\noindent {\bf Abstract.} {\small
We showed in \cite{iv2} that every connected graph can be realized as 
the cut locus of some point on some riemannian surface $S$. 
Here, criteria for the orientability of $S$ are given, and are applied 
to classify the distinct, orientable, cut locus structures on graphs with four generating cycles.
\\{\bf Math. Subj. Classification (2000):} | 53C22, 05C10}


\section{Introduction}

By a surface we always mean a complete, compact and connected $2$-dimen\-si\-o\-nal riemannian manifold without boundary.
Unless explicitly stated otherwise, all our graphs will be connected, undirected, 
and may have multiple edges and loops but not vertices of degree two.

The notion of cut locus was introduced by H. Poincar\'e \cite{p} in 1905, and gain since then an important place in global riemannian geometry.
The {\it cut locus $C(x)$ of} the point $x$ in the riemannian manifold $M$ is the set of all extremities (different from $x$) of maximal (with respect to inclusion) shortest paths starting at $x$; for basic properties and equivalent definitions refer, for example, to \cite{ko} or \cite{sa}.

S. B. Myers \cite{M} established that the cut locus of a real analytic riemannian surface is (homeomorphic to) a graph, 
and a partial converse to his result was proven in \cite{iv2}, namely that every (metric) connected graph can be realized as a cut locus; i.e., there exist a riemannian surface $S_G=(S_G,h)$ and a point $x \in S_G$ such that $C(x)$ is isometric to $G$.
If moreover $G$ is cyclic of order $k$, then it can be realized on a surface of constant curvature \cite{iv2}.

In \cite{iv1} we introduced the notion of cut locus structure on a graph, and discussed its basic properties, while in \cite{iv3} we proposed upper bounds on the number of such structures.
The generic behaviour of cut locus structures is also presented in \cite{iv2}.

In this paper, we are concerned about the orientability of the surfaces $S_G$ realizing the graph $G$ as a cut locus.
If $G$ has an odd number of generating cycles then any surface realizing $G$ is non-orientable.
If the number of generating cycles is even then one cannot generally distinguish, by simply looking to the graph $G$, whether $S_G$ is orientable or not.
In other words, seen as a graph, {\sl the cut locus does not encode the orientablity of the ambient space}; this is our main motivation to introduce ``cut locus structures'' on graphs.
In order to characterize the orientability of a surface $S_G$ realizing $G$ as a cut locus, 
we codify in a cut locus structure any small neighbourhood of the cyclic part $G$ of $C(x)$ in $S_G$.

In Section \ref{CL-str} we briefly present the notion of cut locus structure, 
in Section \ref{crit_orient} we characterize those cut locus structures living on orientable surfaces, while
and in Section \ref{or_graphs} we give several criteria to provide {\it orientably realizable} graphs (i.e., graphs having at least one orientable realization). 
For example, we obtain (Corollary \ref{Petersen}) that the Petersen graph has at least two orientable realizations.

As an application of our results in $\S$\ref{crit_orient}, in Section \ref{S_4cy} we classify all
distinct orientable cut locus structures on graphs with four generating cycles.
In particular, this enables us to point out such a graph which is not orientably realizable.


\bigskip

At the end of this section we recall a few definitions and facts about graphs, to fix notation.

We shall denote by $m(G)$ the number of edges in the graph $G$, and by
$n(G)$ the number of its vertices.

An edge in $G$ is called {\it external} if it is incident to (least) one vertex of degree one, and is called a {\it bridge} if its removal disconnects $G$.

Denote by $B$ the set of all {\it bridges} in the graph $G$. 
Each non-vertex component of $G \setminus B$ is called a $2$-{\it connected component} of $G$.

A $k${\it -graph} is a graph all of whose vertices have degree $k$.
In particular, a $3$-graph will also be called {\it cubic}.
An {\it edge contraction} in the graph $G$ is an operation which removes an edge from $G$ while simultaneously merging together the two vertices it used to connect to a new vertex; all other edges incident to either of the two vertices become incident to the new vertex.

Consider the graph $G$ as a simplicial complex;
the {\it cyclic part} of $G$ is the minimal (with respect to inclusion) subgraph $G^{cp}$ of $G$, 
to which $G$ is contractible; i.e., the minimal subgraph of $G$ obtained by repeatedly contracting external edges, 
and for each vertex remaining of degree two (if any) merging its incident edges.
A graph is called {\it cyclic} if it is equal to its cyclic part.

The power set ${\cal E}$ of $E$ becomes a $Z_2$-vector space over the two-element field $Z_2$, 
if endowed with the symmetric difference $\ast$ as addition, and it is called the {\it edge space} of $G$. 
The {\it cycle space} is the subspace ${\cal Q}$ of ${\cal E}$ generated by (the edge sets of) all the simple cycles of $G$.
If $G$ is seen as a simplicial complex, ${\cal Q}$ is the space of $1$-cycles of $G$ with mod $2$ coefficients.
The symmetric difference $\ast$ of two simple cycles is either a simple cycle or a union of edge-disjoint simple cycles. 
The dimension $q=q(G)$ of the cycle space of the graph $G$ is given by $q(G)=m(G)-n(G)+1$.


\section{CL-structures}
\label{CL-str}

In this section we briefly present the notion of cut locus structure, that we introduced in \cite{iv1}.

\begin{df}
A $G$-{\sl patch} on the graph $G$ is a topological surface $P_G$ with boundary, containing (a graph isometric to) $G$ and contractible to $G$.

A $G$-{\sl strip} is a $G$-patch whose boundary is topologically a circle.

A {\sl cut locus structure} (shortly, a {\sl CL-structure}) on the graph $G$ is a strip on the cyclic part $G^{cp}$ of $G$.

A CL-structure on $G$ is {\sl orientable} if the corresponding $G^{cp}$-strip is an orientable surface.
\end{df}

\begin{df}
An {\sl elementary decomposition} of a $G$-patch $P_G$ is a decomposition of $P_G$ into {\sl elementary strips} such that:
\\- each edge-strip corresponds to precisely one edge of $G$;
\\- each point-strip corresponds to precisely one vertex of $G$.
\end{df}

Denote by ${\cal P}$ and ${\cal A}$ the set of all elementary point-strips, respectively edge-strips, 
of a CL-structur ${\cal C}$e on the graph $G$.

\begin{df}
Consider an elementary decomposition of the $G$-strip $P_G$ such that each elementary strip has a distinguished face, 
labeled $\bar 0$. The face opposite to the distinguished face will be labeled $\bar 1$. 
Here, $\bar 0$ and $\bar 1$ are the elements of the $2$-element group $(Z_2, \oplus)$. 

To each pair $(v,e) \in V \times E$ consisting of a vertex $v$ and an edge $e$ incident to $v$, 
we associate the $Z_2$-sum $\bar s (v,e)$ of the labels of the elementary strips $\nu \in {\cal P}$, 
$\varepsilon \in {\cal A}$ associated to $v$ and $e$; i.e.,
$\bar s (v,e) = \bar 0$ if the distinguished faces of $\nu$ and $\varepsilon$ agree to each other, and $\bar 1$ otherwise.
Therefore, to any cut locus structure ${\cal C}$ we can associate a function $s_{\cal C} :E \to \{\bar 0,\bar 1\}$, 
\begin{eqnarray}
\label{companion}
s_{\cal C}(e)= \bar s (v,e) \oplus \bar s (v',e),
\end{eqnarray}
where $v$ and $v'$ are the vertices of the edge $e \in E$.

We call the function $s_{\cal C}$ defined by (\ref{companion}) the {\sl companion function} of ${\cal C}$. 

An edge-strip $P_e$ (or simply an edge $e$) in a CL-structure ${\cal C}$ is called {\sl switched} if $s_{\cal C} (e)=\bar 1$.
\end{df}

\begin{df}
\label{diff_CL}
The {\sl CL-structures} ${\cal C}$, ${\cal C}'$ on $G$ are called {\sl equivalent} 
if their characteristic functions are equivalent on every $2$-connected component $H$ of $G$; 
i.e., on every $H$ either $s_{\cal C} = s_{{\cal C}'}$, or $s_{\cal C} =\bar 1 \oplus s_{{\cal C}'}$.
\end{df}

We shall use the following way to planary represent a CL-structure ${\cal C}$ on the $3$-graph $G$ \cite{iv1}:
we represent in the plane each vertex-strip such that its distinguished face is ``up'', 
and afterward connects the vertex-strips by edge-strips.
To schematically represent this, we shall overwrite an ``{\rm x}'' to the drawn image of an edge if its strip is switched, 
and an ``{\rm =}'' to the drawn image of a edge if its strip is not-switched.
If, moreover, $G$ is planar then one can use any planar representation of $G$
to schematically represent any $G$-strip.

\bigskip

We explain now the relationship between patches and cut locus realizations of graphs.

Assume that the cut locus $C(x)$ of the point $x$ in the surface $S$ is isometric to the graph $G$.
Then, cutting off the surface an open intrinsic disc of radius smaller than the injectivity radius at $x$,
provides a strip on $G$, and consequently on $G^{cp}$, called 
{\it the cut locus natural structure of} $x$, and denoted by $CLNS(x)$. 

The converse is established by the following result, allowing us to consider strips whenever we think about cut locus realizations of graphs.

\begin{thm}
\label{glue} {\rm \cite{iv2}}
For every graph $G$ there exists at least one $G$-strip, and 
each $G$-strip provides a realization of $G$ as a cut locus.
\end{thm}

We shall implicitly use the following simple result, easily obtained from the above considerations.

\begin{lm}
\label{cy-sw}
Let $G$ be a $3$-graph. In any planar representation of a $G$-strip, 
each cycle-patch contains at least one switched edge-strip.
\end{lm}


\section{Orientable cut locus structures}
\label{crit_orient}

We are concerned next about the orientability of the surfaces $S_G$ realizing the graph $G$ as a cut locus.

If $G$ has an odd number of generating cycles then any surface realizing $G$ is non-orientable, 
because any cut locus on an orientable surface of genus $g$ has $2g$ generating cycles.

If the number of generating cycles is even then one cannot generally distinguish, 
by simply looking to $G$, whether $S_G$ is orientable or not, see Example \ref{Non_orient} or Theorem \ref{4}.

\bigskip

The following result expresses formally the simple fact that a circle-patch is an orientable surface if and only if
it has an even number of switches.
Together with Theorem \ref{or_cy}, it will be repeatedly applied in Section \ref{S_4cy}.

\begin{lm}
\label{sum_e}
A patch $P$ over a cycle $C$ is an orientable surface if and only if 
\begin{equation}
\label{oplus_cycle_patch}
\oplus_{e \in E(C)} s(e)= \bar 0 .
\end{equation}
\end{lm}

\begin{lm}
\label{prod_or_cy}
The product $\ast$ of cycles defines a natural operation for the cycle-patches.
In particular, if the cycle $C$ is the product of the cycles $C_1, ..., C_k$, 
$C=C_1 \ast ...\ast C_k$, each of which given with an orientable patch, then the induced $C$-patch is orientable.
\end{lm}

\noindent{\sl Proof:}
We prove the result by induction over $k$.

Assume first $k=2$, and let $\varepsilon_1$ be the sum of switches in $C_1 \setminus C_2$, $\varepsilon_2$ be the sum of switches in $C_2 \setminus C_1$, and $\varepsilon_{12}$ the sum of switches in $C_1 \cap C_2$.

By Lemma \ref{sum_e}, the $C_1$-patch is orientable if and only if
\begin{equation}
\label{C1}
\varepsilon_1 \oplus \varepsilon_{12} = \bar 0,
\end{equation}
the $C_2$-patch is orientable if and only if
\begin{equation}
\label{C2}
\varepsilon_2 \oplus \varepsilon_{12} = \bar 0,
\end{equation}
and respectively the $C=C_1 \ast C_2$-patch is orientable if and only if
\begin{equation}
\label{C1C2}
\varepsilon_1 \oplus \varepsilon_2 = \bar 0.
\end{equation}
But the equation (\ref{C1C2}) follows by simply adding (\ref{C1}) and (\ref{C2}).

Assume now that the statement is true for $k$ cycle patches; in order to prove it for $k+1$ cycle patches, simply put
$C'_1=C_1 \ast...\ast C_k$, $C'_2=C_{k+1}$, and apply the case $k=2$ to $C'_1$ and $C'_2$.
\hfill $\Box$

\bigskip

Recall that for any surface $S$ and any point $x$ in $S$, $C(x)$, if not a single point, is a local tree (i.e., each of its points $z$ has a neighbourhood $V$ in $S$ such that the component $K_z(V)$ of $z$ in $C(x)\cap V$ is a tree).
A {\it tree} is a set $T$ any two points of which can be joined
by a unique Jordan arc included in $T$.

Eventhough a cut locus $C(x)$ may be quite large a set
(see \cite{GS} or \cite{H2} for examples of non triangulable cut loci),
its {\it cyclic part} is a cyclic graph, see for example \cite{Itoh-Zamfirescu}.

The {\it tangential cut locus} of the point $x \in S$ is the boundary of
the maximal (with respect to inclusion) domain in the tangent space $T_xS$ to $S$ at $x$, to which the exponential map at $x$ is a diffeomorphism.

The last part of the following preliminary result has some interest in its own.

\begin{lm}
\label{n-or_crit}
A surface $S$ is non-orientable if and only if for any point $x$ in $S$ there exists an edge $e$ of $C(x)^{cp}$ whose two images in the tangential cut locus have the same orientation. Moreover, such an edge $e$ cannot be a bridge of $C(x)^{cp}$.
\end{lm}

\noindent{\sl Proof:}
The ``if and only if'' part is clear.

For the last part, fix some point $x$ in $S$ and assume $e=12$ is a bridge in $G=C(x)^{cp}$, the removal of which yields two subgraphs $G_1$, $G_2$ of $G$, with $i$ a vertex of $G_i$, $i=1,2$.
Let $P_i$ be the patch induced by $CLNS(x)$ along $G_i$, $i=1,2$.

Now regard the boundary of $CLNS(x)$ as a directed curve $O$ homeomorphic to a circle. 
$O$ enters $P_1$ at $1$, covers twice all edges in $G_1$, and
exits $P_1$ again at $1$; afterwards it goes along $e$ to $2$, 
enters $P_2$ at $2$, covers twice all edges in $G_2$, and
exits $P_2$ again at $2$. Therefore, with $O_i=O \cap P_i$, $i=1,2$, 
the order along $O$ is $e=21$, $1$, $O_1$, $1$, $e=12$, $2$, $O_2$, $2$.
Since this is also the corresponding order in the tangential cut locus, the proof is complete.
\hfill $\Box$

\bigskip

One can roughly figure out the following statement by the fact that a non-orientable surface has only one ``face''.

\begin{thm}
\label{crit}
The surface $S_G$ realizing the connected graph $G$ as a cut locus of the point $x$ in $S$ is non-orientable if and only if there exists a cycle $C$ of $G$ and a non-orientable $C$-patch induced by $CLNS(x)$.
\end{thm}

\noindent{\sl Proof:}
By Lemma \ref{n-or_crit}, the surface $S_G$ is non-orientable if and only if there exists an edge $e$ in $G$ whose two images in the tangential cut locus of $x$ have the same orientation.
Again by Lemma \ref{n-or_crit}, such an edge can always be included in a cycle, whose induced patch is consequently non-orientable.
\hfill $\Box$

\begin{co}
\label{loop_at_3}
If the cyclic part of $G$ has a loop at a degree three vertex then $S_G$ is non-orientable.
\end{co}

\noindent{\sl Proof:} All CL-structures on $G$ induce the same patch along a loop at a degree three vertex of $G$,
which is non-orientable by Lemma \ref{sum_e}, and therefore $S_G$ is non-orientable by Theorem \ref{crit}.
\hfill $\Box$

\begin{ex}
\label{loop_at_4}
There exist orientable realizations of graphs having loops at vertices of degree four (or more).
To see this, consider a flat torus $T$ of square fundamental domain $D$. 
Denote by $x$ the point corresponding to the vertices of $D$, and consider $CLNS(x)$.
The cyclic part of $C(x)$ consists of two loops at a degree four vertex and $T$ is orientable, 
see Figure 1.
\end{ex}

\begin{figure*}
\label{Orient_Cy_deg_4}
\centering
  \includegraphics[width=0.5\textwidth]{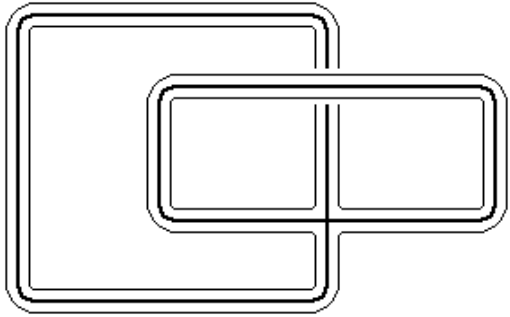}
\caption{CLNS on a flat torus of square fundamental domain.}
\end{figure*}

Theorem \ref{crit} can alternatively be stated as follows.

\begin{thm}
\label{or_cy}
The surface $S_G$ realizing the connected graph $G$ as a cut locus is orientable if and only if there exists a system of generating cycles for $G^{cp}$, each of which has an orientable patch in $S_G$.
\end{thm}

\noindent{\sl Proof:}
If there exist an orientable surface $S_G$ realizing $G$ as a cut locus then
every cycle of $G$ has an orientable patch in $S_G$, by Theorem \ref{crit}.

Conversely, assume there exists a system of generating cycles for $G^{cp}$, each of which has an orientable patch in the surface $S_G$ realizing $G$ as a cut locus. Then, by Lemma \ref{prod_or_cy}, each cycle in $S_G$ is orientable.
\hfill $\Box$

\begin{ex}
\label{Non_orient}
By the use of Theorems \ref{crit} and \ref{or_cy} one can easily see that 
the CL-structure in Figure 2 (a) is orientable, while the one in Figure 2 (b) is not.
To extend this example to a $3$-graph with $2g$ generating cycles, 
simply connect by an edge-strip the non-orientable CL-structure in Figure 2 (b) 
to the orientable CL-structure of the graph with $2g-2$ in Figure 2 (c), see Figure 2 (d).
\end{ex}

\begin{figure*}
\label{No_orient}
\centering
  \includegraphics[width=0.6\textwidth]{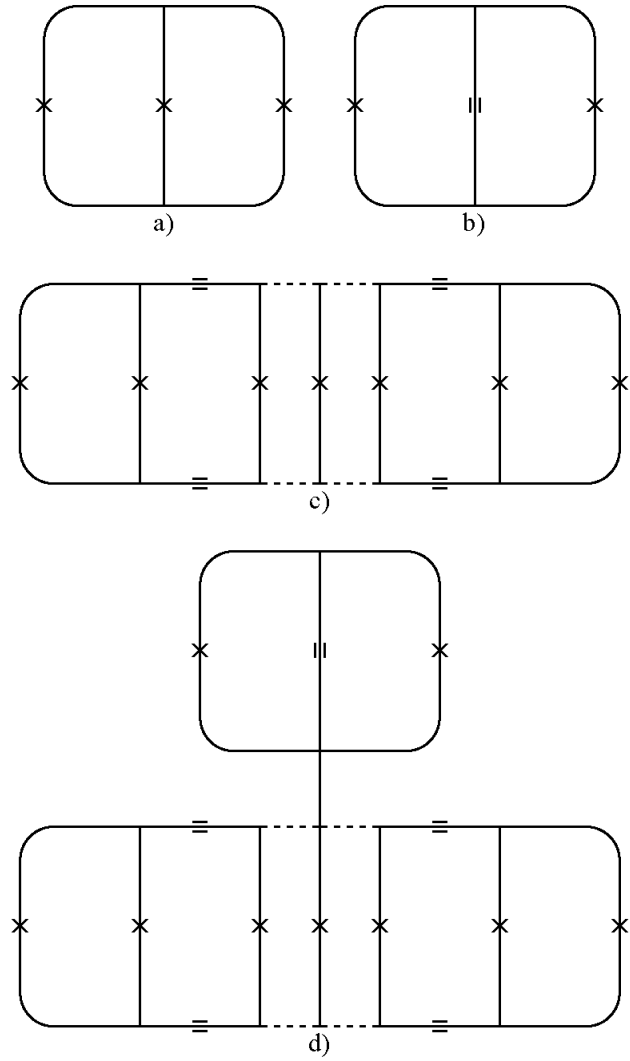}
\caption{Cut locus structures: orientable (a) and (c), and non-orientable (b) and (d).}
\end{figure*}

The following simple remark is called a corollary because of it importance for the next section.

\begin{co}
\label{2v}
Let $S_G$ be an orientable realization of the $3$-graph $G$, reprezen\-ted planary.
Then, for every cycle $C$ of $G$, if $C$ consists of two edges then both of them are switched, 
and if $C$ consists of three edges then exactly two of them are switched.
\end{co}


\section{Orientable realizations of small graphs}
\label{S_4cy}

As an application of our previous results, we present in this section all 
distinct, orientable cut locus structures on $3$-graphs with four generating cycles;
the others CL-structures realized as cut loci on surfaces of genus $2$ 
can be obtained from those on $3$-graphs with $4$ generating cycles by contracting edge-strips, see \cite{iv1}.

The following statement can be obtained by straightforward inductive constructions.

\begin{figure*}
\centering
  \includegraphics[width=0.8\textwidth]{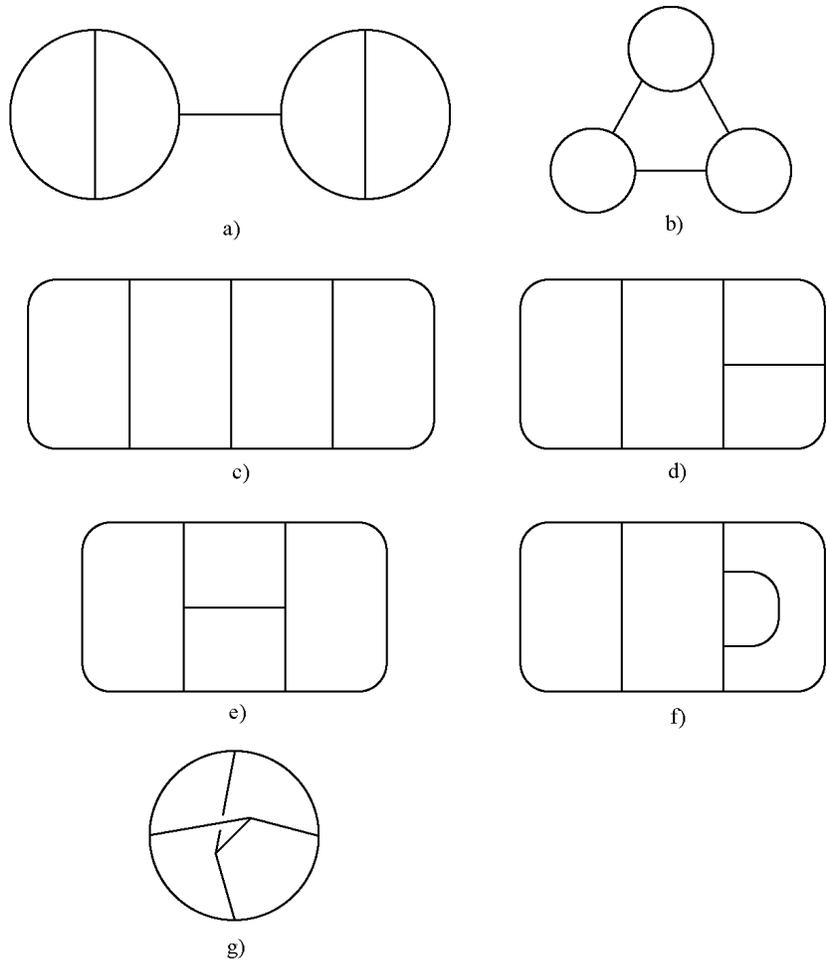}
\caption{All $3$-graphs with $4$ generating cycles and no loops.}
\label{4cy}
\end{figure*}

\begin{lm}
\label{g43}
There are only $7$ (up to isomorphisms) $3$-graphs with $4$ generating cycles and no loops; 
they are listed in the Figure \ref{4cy}.
\end{lm}

\begin{thm}
\label{4}
The list of all orientable cut locus structures on $3$-graphs with $4$ generating cycles 
is presented in the Figures 4 -- 13.
\end{thm}

\noindent{\sl Proof:}
By Theorem \ref{or_cy} and Lemma \ref{sum_e}, it suffices to label each edge of our graphs with $\bar 0$ and $\bar 1$ in such a manner that for each simple generating cycle $C$, $\oplus_{e \in E(C)} s(e)= \bar 0$.

The proofs for the planar graphs in Figure \ref{4cy} (a)-(f) are similar, 
as they all employ the planar representation of any CL-structure on these graphs (see Lemma \ref{cy-sw}).
For such structures, by Corollary \ref{2v}, both edges of any cycle consisting of two edges will be switched, 
and precisely two edges of any cycle consisting of three edges will be switched.

\begin{figure*}
\label{4cy-i}
\centering
  \includegraphics[width=0.50\textwidth]{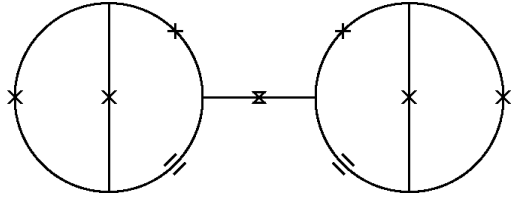}
\caption{Unique orientable cut locus structure for the graph in Figure \ref{4cy} (a).}
\end{figure*}

The graph $G_a$ in Figure \ref{4cy} (a) is symmetric with respect to (the mid-point of) its bridge.
Both cycles of two edges have their edges switched, while the two cycles of three edges have one edge not-switched.
Since the bridge cannot modify a CL-structure (Definition \ref{diff_CL}), 
we obtain the unique result illustrated in Figure \ref{4cy-i}.

\begin{figure*}
\label{4cy-ii}
\centering
  \includegraphics[width=0.4\textwidth]{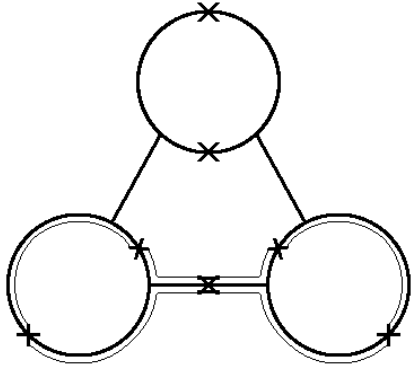}
\caption{No orientable cut locus structure for the graph in Figure \ref{4cy} (b).}
\end{figure*}

Assume the graph $G_b$ in Figure \ref{4cy} (b) has at least one orientable CL-structure. Then, all edges of the three cycles of two edges in $G_b$ are labeled $\bar 1$. But these labels contradict the existence of a $G_b$-strip, because they force a $G_b$-patch to have three boundary components, see Figure \ref{4cy-ii}.

\begin{figure*}
\label{4cy-iii}
\centering
  \includegraphics[width=1.0\textwidth]{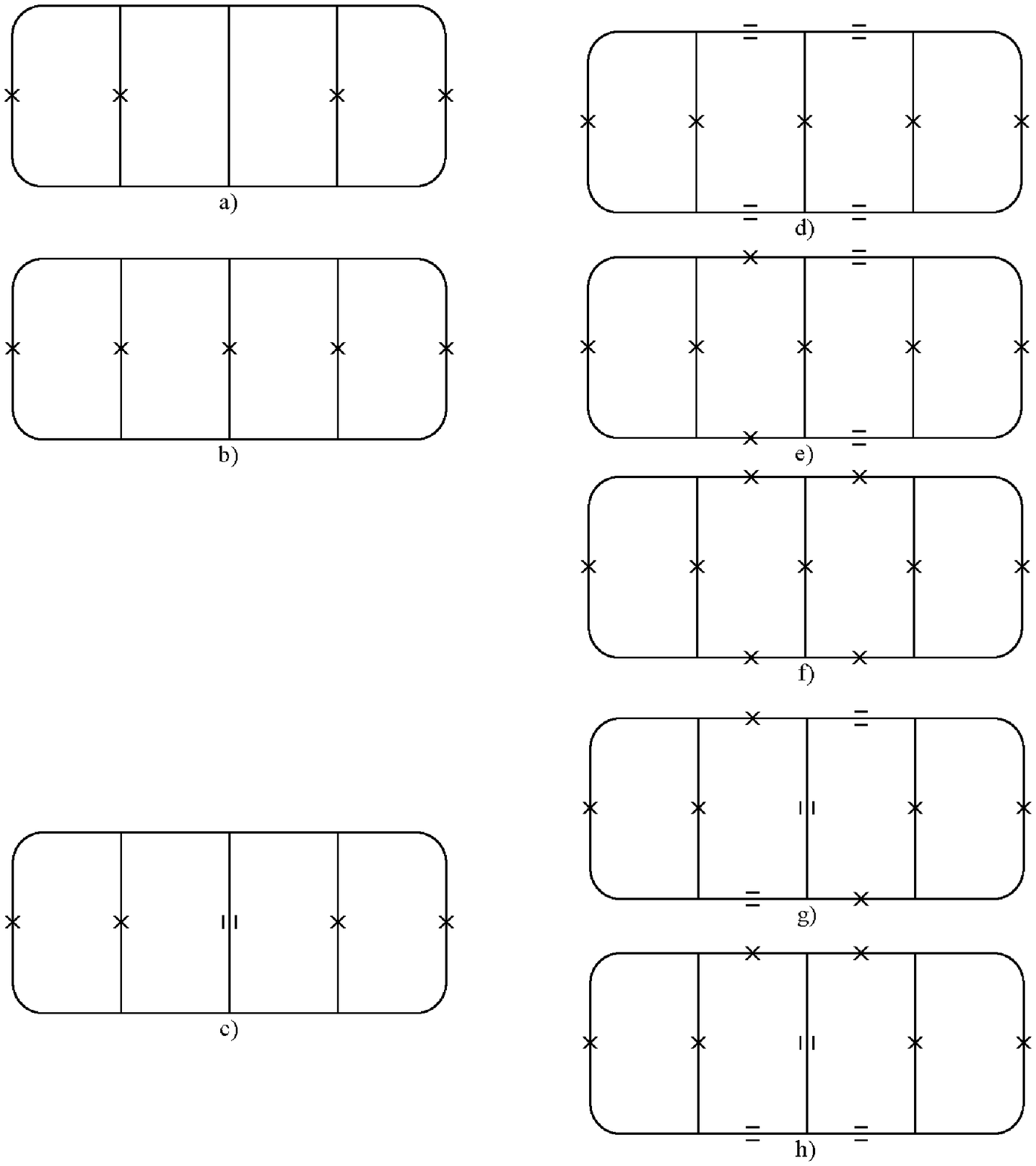}
\caption{Five orientable cut locus structures for the graph in Figure \ref{4cy} (c).}
\end{figure*}

For the graph $G_c$ in Figure \ref{4cy} (c), we first label by $\bar 1$ all edges of its two cycles of two edges, see Figure \ref{4cy-iii} (a). Figure \ref{4cy-iii} (b)-(c) shows the next step of our labeling.
Starting from Figure \ref{4cy-iii} (b), we have three possibilities to label the remaining edges according to Lemma \ref{sum_e}, shown in Figure \ref{4cy-iii} (d)-(e)-(f).
Starting from Figure \ref{4cy-iii} (c) and taking into account the symmetries of $G_c$, we have another two possibilities to label the remaining edges according, shown in Figure \ref{4cy-iii} (g)-(h).
By Theorem \ref{or_cy} and Lemma \ref{sum_e}, all obtained CL-structures (Figure \ref{4cy-iii} (d)-(h)) are orientable.

\begin{figure*}
\label{4cy-iv}
\centering
  \includegraphics[width=1.0\textwidth]{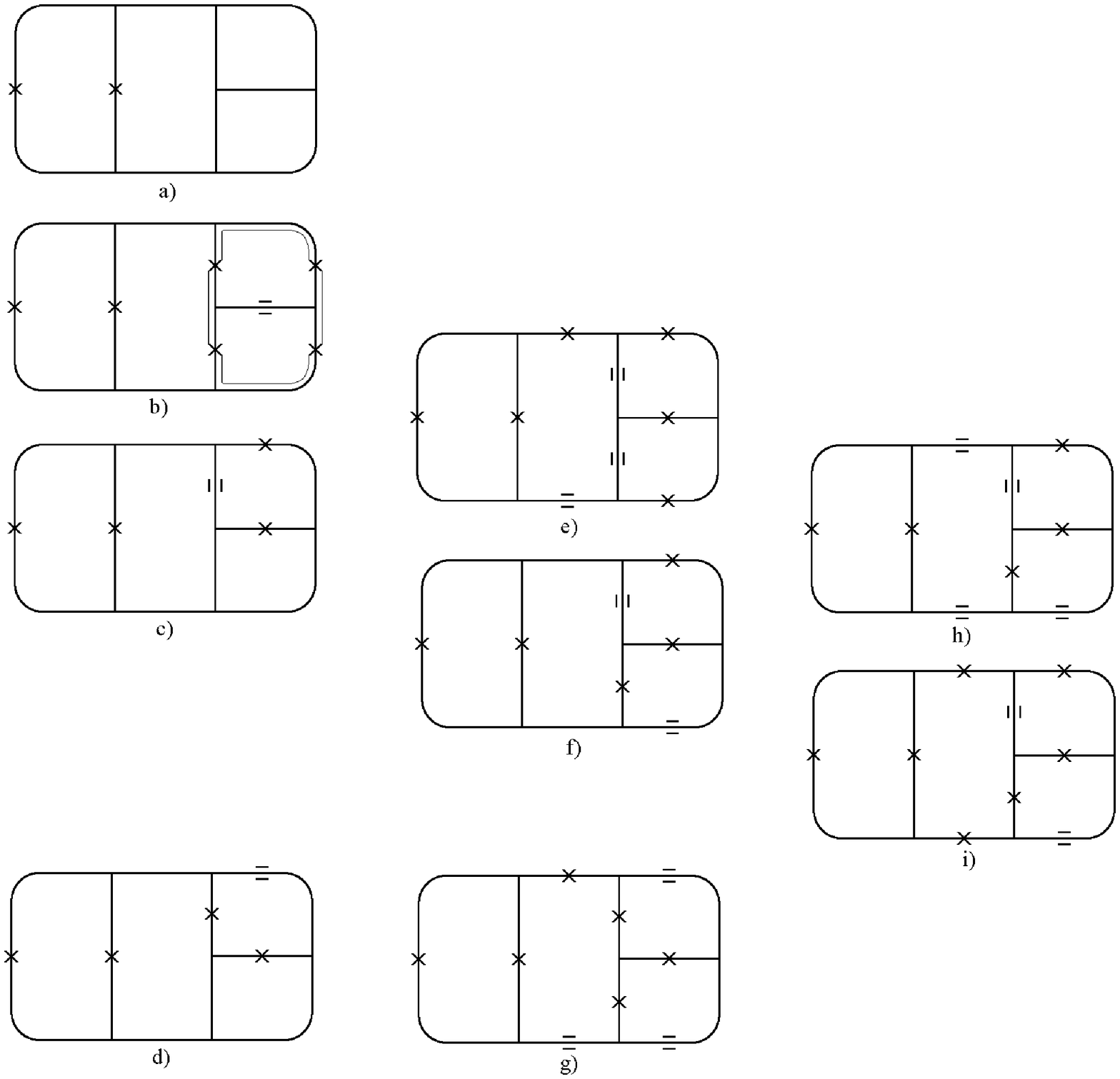}
\caption{Four orientable cut locus structures for the graph in Figure \ref{4cy} (d).}
\end{figure*}

For the graph $G_d$ in Figure \ref{4cy} (d), we apply first Corollary \ref{2v} to obtain Figure \ref{4cy-iv} (a). Then we consider the cases in Figure \ref{4cy-iv} (b), (c), and (d), according to the labeling of cycles with three edges.
The case in Figure \ref{4cy-iv} (b) provides a patch which is not a strip, impossible.
For Figure \ref{4cy-iv} (c) we have two subcases to treat, illustrated in Figure \ref{4cy-iv} (e) and (f), and (f) further ramificates to (h) and (i).
Excluding the subcases of (c), (d) produces one more CL-structure, see Figure \ref{4cy-iv} (g).

\begin{figure*}
\label{4cy-v}
\centering
  \includegraphics[width=0.62\textwidth]{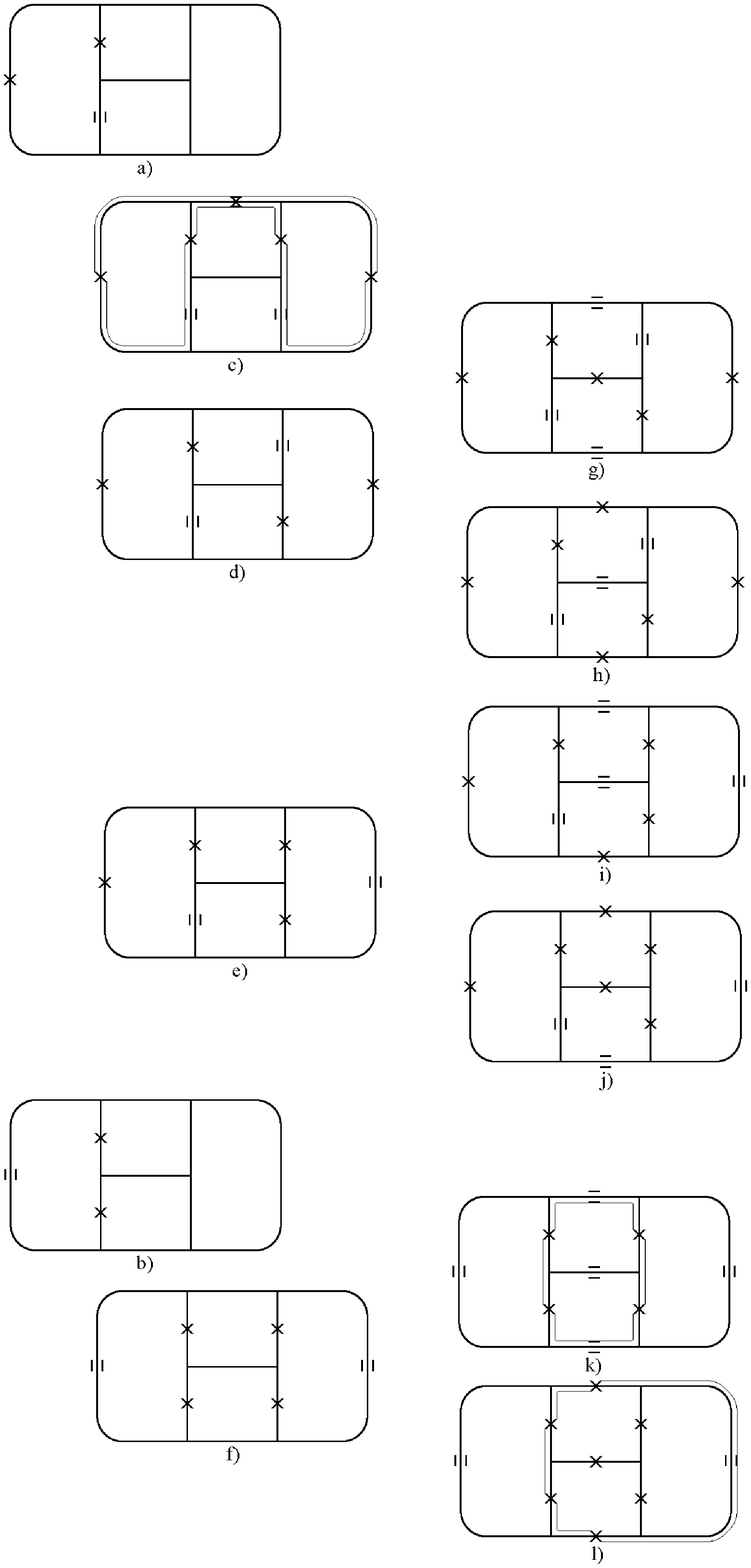}
\caption{Four orientable cut locus structures for the graph in Figure \ref{4cy} (e).}
\end{figure*}

The symmetries of the graph $G_e$ in Figure \ref{4cy} (e) help to reduce the number of considered subcases.
At the beginning we have to take into account two cases, illustrated in Figure \ref{4cy-v} (a) and (b), 
according to the edges of the first cycle of three edges.
Labeling the second cycle of three edges in Figure \ref{4cy-v} (a) produces the subcases (c), (d), and (e), 
of which (c) is not a strip, while (d) and (e) provide each two CL-structures (Figure \ref{4cy-v} (g)-(h) and (i)-(j)).
Excluding the previously treated subcase (e), from (b) we obtain (f) and further (k) and (l), both of which are not strips.

\begin{figure*}
\label{4cy-vi}
\centering
  \includegraphics[width=1.0\textwidth]{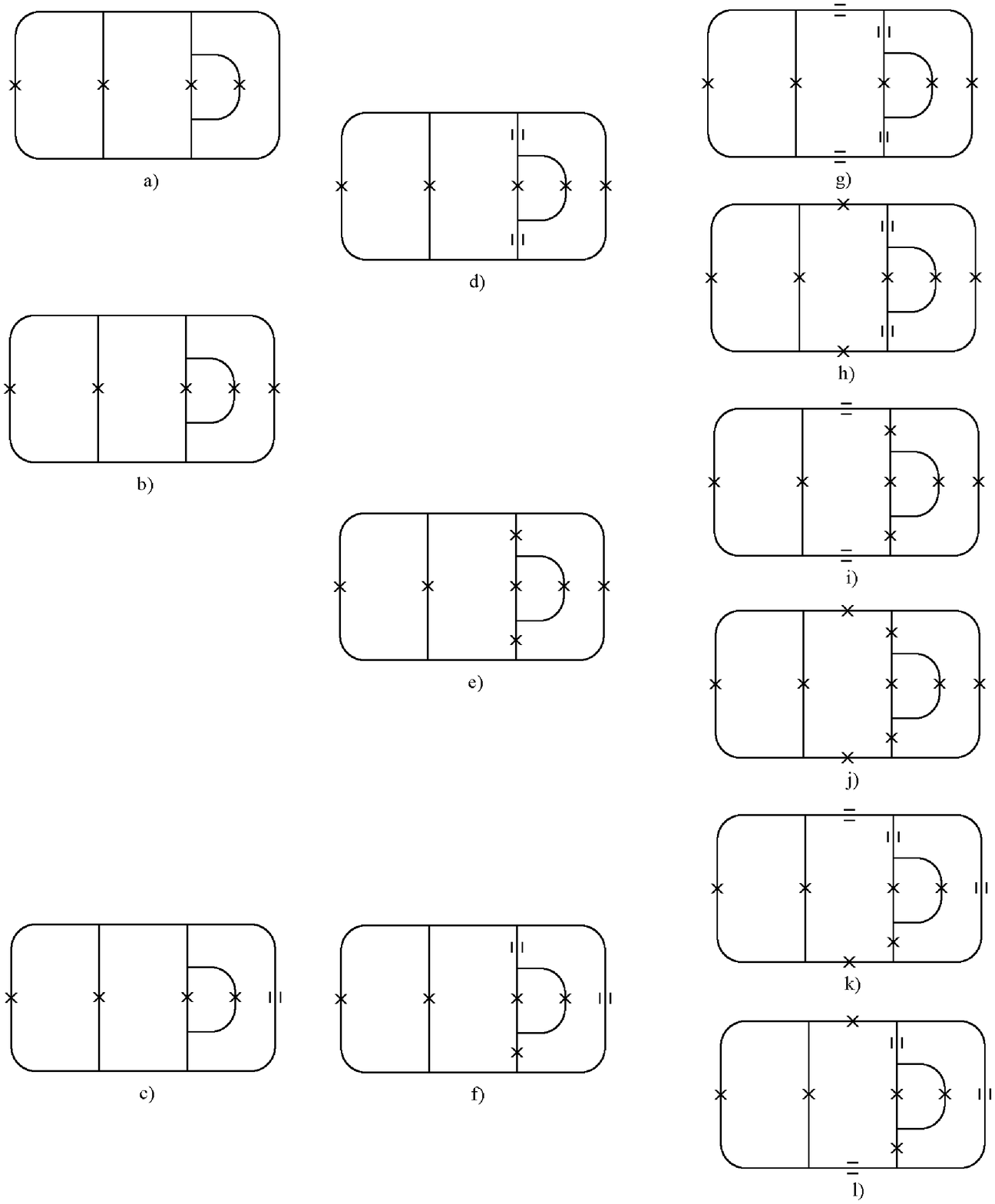}
\caption{Six orientable cut locus structures for the graph in Figure \ref{4cy} (f).}
\end{figure*}

The graph $G_f$ in Figure \ref{4cy} (f) has two cycles of two edges, see Figure \ref{4cy-vi} (a).
There we have to treat the two subcases (b) and (c), (b) yielding (d) and (e), and (c) yielding (f).
Each of (d), (e), (f) produces two CL-structures, see Figure \ref{4cy} (g)-(h), (i)-(j), and (k)-(l).

\begin{figure*}
\label{4cy-vii_1}
\centering
  \includegraphics[width=0.84\textwidth]{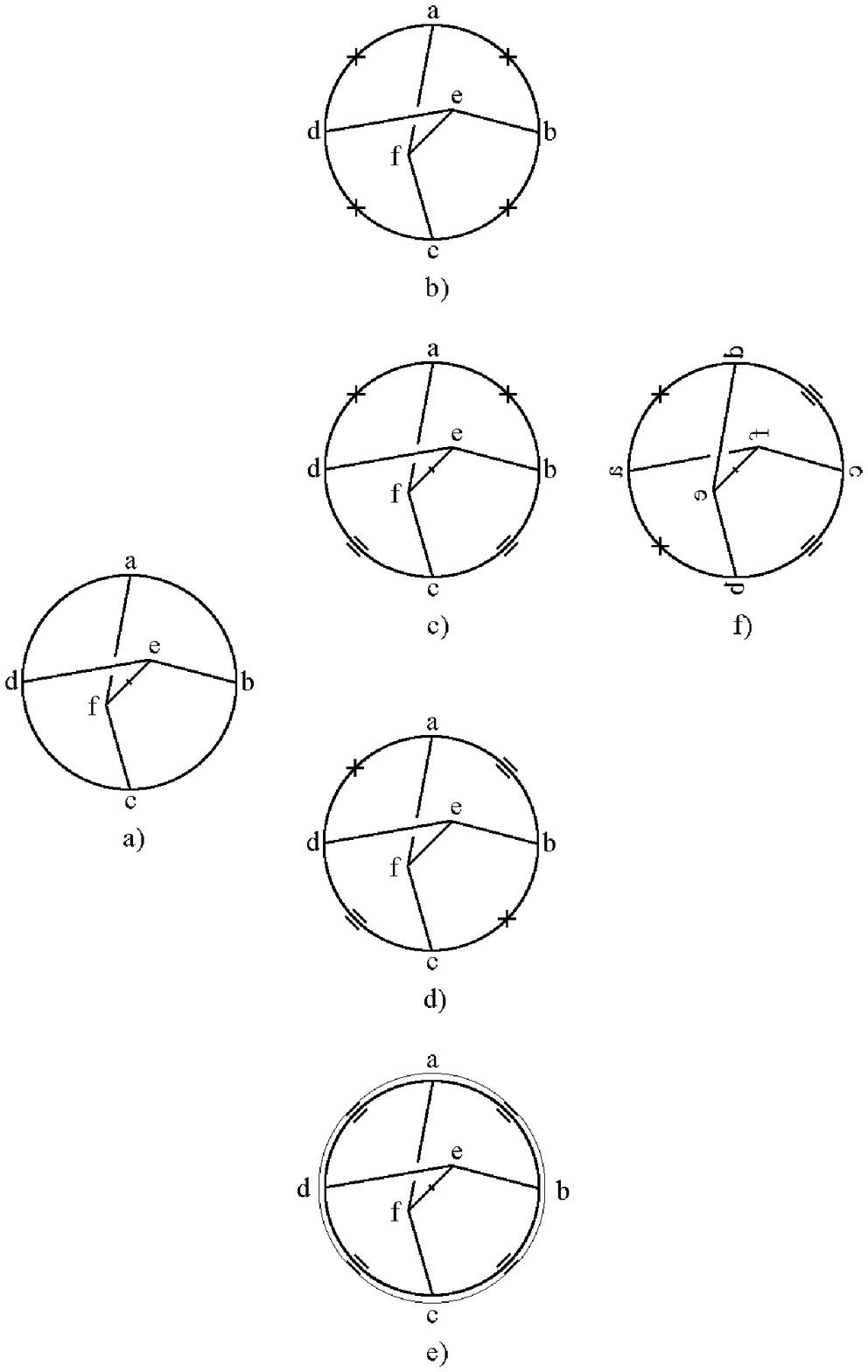}
\caption{Finding orientable cut locus structures for the graph in Figure \ref{4cy} (g); $3$ cases to further consider: (b), (c), and (d).}
\end{figure*}

The graph $G_g$ in Figure \ref{4cy} (g) is not planar, and the analysis  is a little more complicated.
Lebeling the four edge cycle $abcd$ on the outer circle in  Figure \ref{4cy-vii_1} (a) produces
five cases. Two of them, (c) and (f) in Figure \ref{4cy-vii_1}, are identical up to the central symmetry of $G_g$ with respect to the mid-point of the edge $ef$, and one other ((e) in Figure \ref{4cy-vii_1}) more 
than one boundary component, impossible, hence there remain only (b), (c) and (d) in Figure \ref{4cy-vii_1}.

\begin{figure*}
\label{4cy-vii_2}
\centering
  \includegraphics[width=0.95\textwidth]{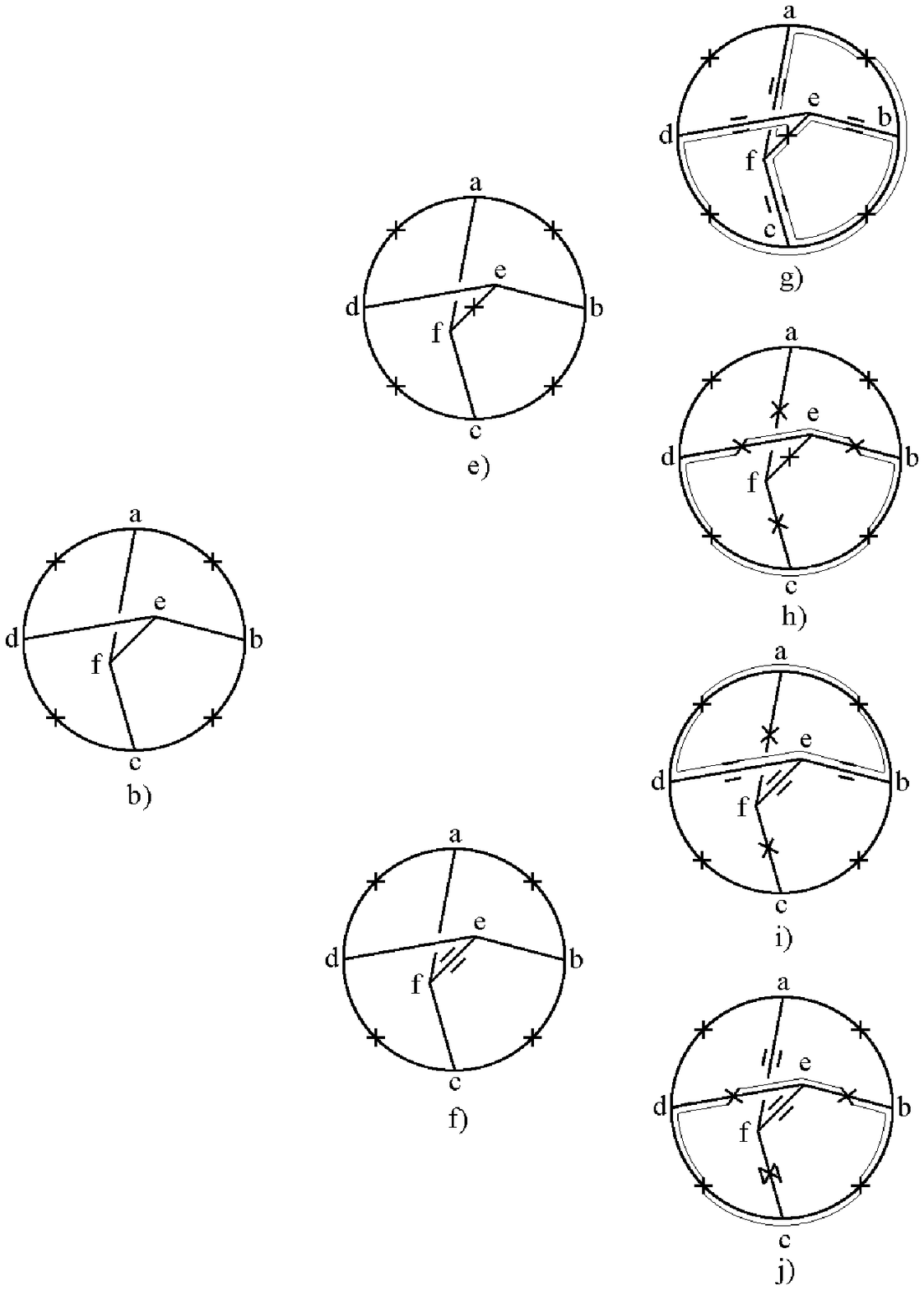}
\caption{No orientable cut locus structure for Figure \ref{4cy-vii_1} (b).}
\end{figure*}

The subcases obtained from the case (b) are illustrated in Figure \ref{4cy-vii_2},
all of them leading to several boundary components.

\begin{figure*}
\label{4cy-vii_3}
\centering
  \includegraphics[width=1.0\textwidth]{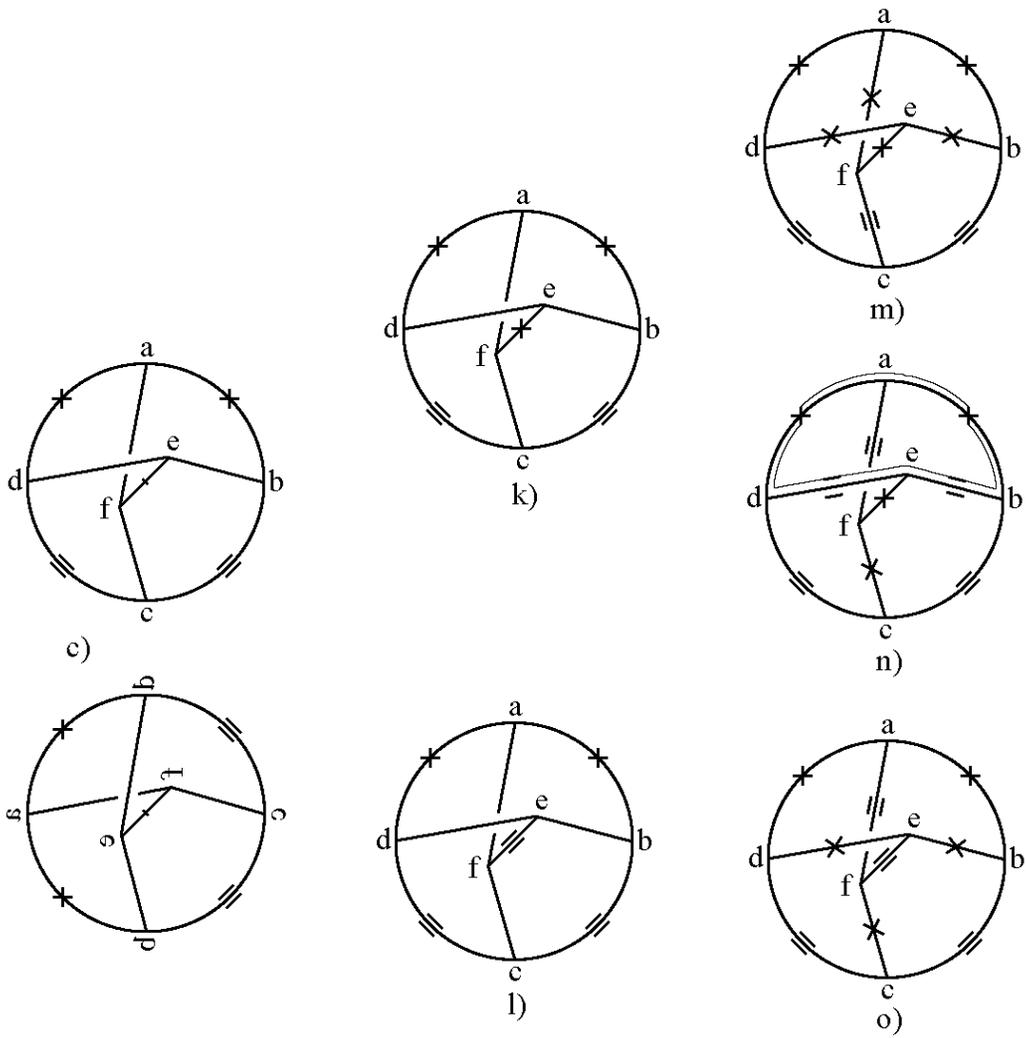}
\caption{Two orientable cut locus structures for Figure \ref{4cy-vii_1} (c).}
\end{figure*}

The subcases obtained from the case (c) are illustrated in Figure \ref{4cy-vii_3},
the first two of them, (m) and (n), leading to several boundary components. The only orientable CL-structure
obtained in this case is represented in Figure \ref{4cy-vii_3} (o).

\begin{figure*}
\label{4cy-vii_4}
\centering
  \includegraphics[width=1.0\textwidth]{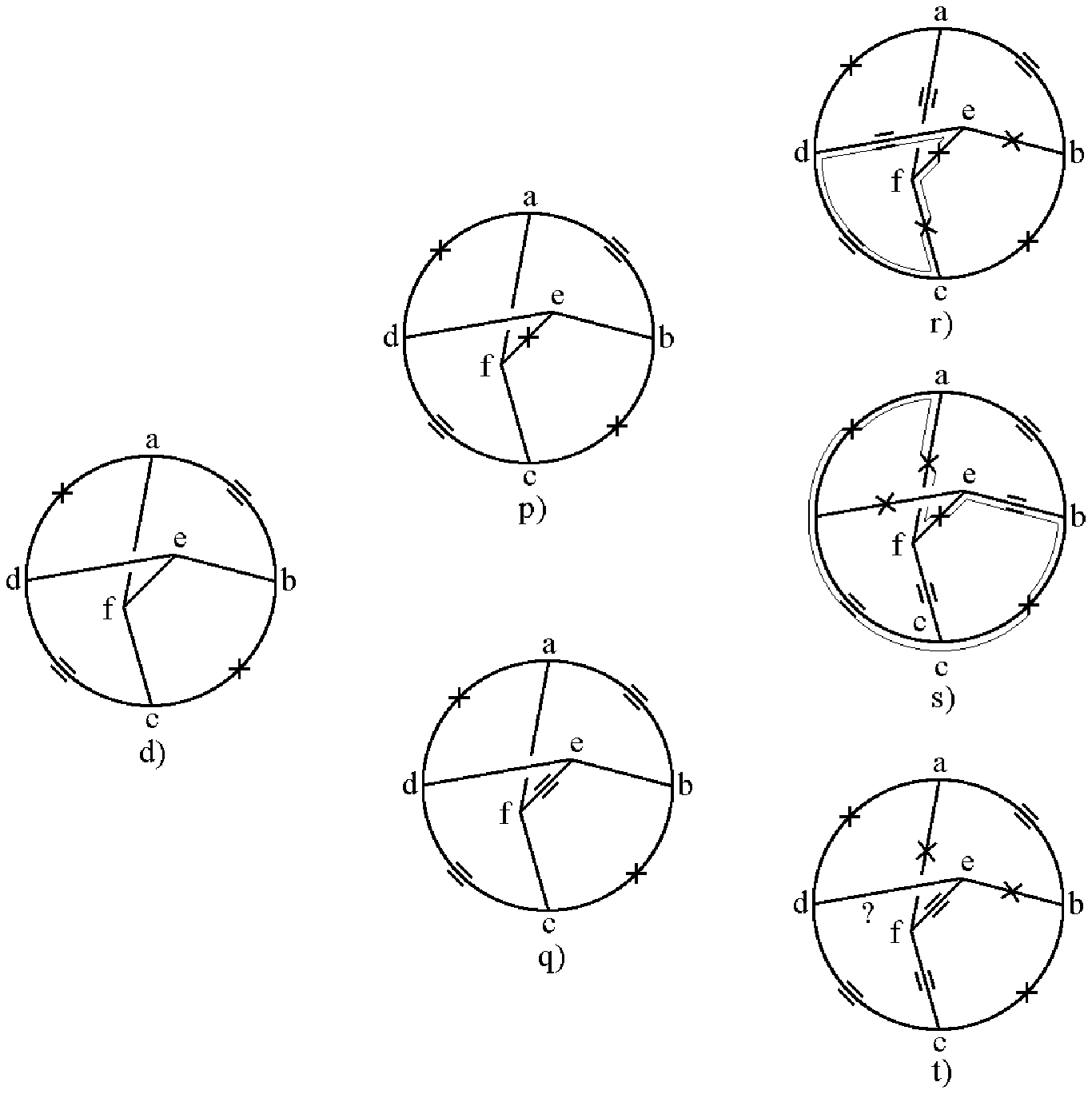}
\caption{No orientable cut locus structure for Figure \ref{4cy-vii_1} (d).}
\end{figure*}

The subcases obtained from the case (d) are illustrated in Figure \ref{4cy-vii_4},
the first two of them, (r) and (s), leading to several boundary components.
The subcase (t) also yields a contradiction, producing either several boundary components, 
or a non-orientable cycle-patch, according to the labeling of the edge $de$.
\hfill $\Box$

\bigskip

Consider a CL-structure ${\cal C}$ on the graph $G$, a riemannian surface $(S,g)$ and a point $x \in S$.
${\cal C}$ is called {\sl stable} with respect to $x \in S$ if (i) $CLNS(x)={\cal C}$, and (ii) there exists a neighbourhood of $x$ in $S$, for all points $y$ of which holds $CLNS(y)={\cal C}$.
The CL-structure ${\cal C}$ is called {\sl stable} if it is stable on all surfaces where it can be realized as a CLNS \cite{iv2}.

We proved in \cite{iv2} that a CL-structure on the graph $G$ is stable if and only if $G$ is cubic; this and Theorem \ref{4} directly imply the following.

\begin{co}
Up to graph homeomorphisms and CL-structures equi\-va\-lence, there are 22 stable and orientable CL-structures on surfacees of genus $2$.
\end{co}


\section{Orientably realizable graphs}
\label{or_graphs}

We noticed that the orientability of a surface $S_G$ realizing the graph $G$ as a cut locus 
is implied not by the properties of $G$, but by those of the cut locus structure of $G$, see Figure 2.
Nevertheless, some graphs have no orientable realization, while others have several, see Section \ref{S_4cy}. 
We present next some classes of such graphs.

\begin{df}
A graph is called {\sl orientably realizable} if there exists at least one orientable riemannian surface 
realizing it as a cut locus.
\end{df}

A first --obvious-- obstruction for a graph being orientably realizable is provided by the odd number of generating cycles.
Another obstruction is the existence of loops at degree three vertices, see Corollary \ref{loop_at_3}.
Theorem \ref{4} shows, in particular, that these two obstructuctions are not sufficient to characterize
the orientably realizable graphs, see Figure 5.

In this section we give a few criteria for orientability.

\begin{co}
\label{tree}
Assume there exist orientably realizable subgraphs $G_1$,..., $G_m$ in the graph $G$ whose union is $G$, 
such that the intersection of any two of them has at most one point. 
If the induced incidence graph of $\{ G_1,..., G_m \}$ is a tree then $G$ is orientably realizable.
\end{co}

\noindent{\sl Proof:} Theorem \ref{or_cy} and induction over $m$.
\hfill $\Box$

\begin{ex}
Corollary \ref{tree} is not necessarily true if the induced incidence graph of $\{ G_1,..., G_m \}$ in $G$ is not a tree, 
as the example of a triangle ($m=3$, each $G_i$ coincides with an edge) shows.
\end{ex}

\begin{co}
\label{conn}
Assume there exist orientably realizable disjoint subgraphs $G_1$, $G_2$ of the graph $G$ such that 
$G \setminus (G_1 \cup G_2)$ is the disjoint union of $k$ edges, each of which connects $G_1$ to $G_2$.
If $k \in \{1, 3\}$ then $G$ is orientably realizable.
\end{co}

\noindent{\sl Proof:} 
The case $k=1$ follows easily from Corollary \ref{tree}, by considering the connecting edge as a third graph.

For $k=3$, put $\{ e_1, e_2, e_3 \}=G \setminus (G_1 \cup G_2)$.
We may assume that the edges $e_1, e_2, e_3$ determine two generating cycles of $G$, say $C$ containing $e_1, e_2$ and $C'$ containing $e_2, e_3$.
Denote by $\varepsilon$ the (mod 2) sum of the switched edges of $C$ in the $G_1$-strip $P_{G_1}$ and the $G_2$-strip $P_{G_2}$; define similarly $\varepsilon'$ for $C'$.
There are three cases to analize.

{\sl i)} If $\varepsilon = \varepsilon' =0$ then join $P_{G_1}$ to $P_{G_2}$ by $3$ switched edge-strips.

{\sl ii)} If $\varepsilon =0$ and $\varepsilon' =1$ (or vice versa) then join $P_{G_1}$ to $P_{G_2}$ by switched $e_1$ and $e_2$-strips and a non-switched $e_3$-strip.

{\sl iii)} If $\varepsilon = \varepsilon' =1$ then join $P_{G_1}$ to $P_{G_2}$ by switched $e_1$- and $e_3$-strips and a non-switched $e_2$-strip. 

It is straightforward to check that the result is, in all cases, a $G$-strip,
which is orientable by Theorem \ref{or_cy}.
\hfill $\Box$

\begin{ex}
Corollary \ref{conn} is not necessarily true if the subgraphs $G_1, G_2$ are not orientably realizable.
To see this, consider the graph $G$ composed of two cycles $G_1, G_2$ joined by an edge ($k=1$); 
by Lemma \ref{cy-sw} and Corollary \ref{2v}, $G$ is not orientably realizable.
\end{ex}

\begin{figure*}
\centering
  \includegraphics[width=0.95\textwidth]{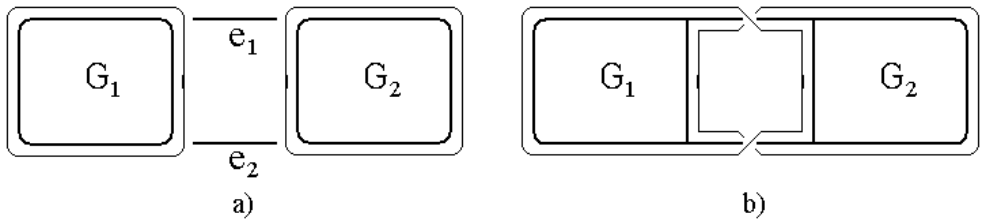}
\caption{Corollary \ref{conn} is not true for $k=2$.}
\label{k=2}
\end{figure*}

\begin{ex}
Corollary \ref{conn} is not necessarily true for $k=2$.
With the notation in its proof,
if $\varepsilon = \varepsilon' =0$ then, imposing to have an orientable cycle along (the edges corresponding to) $\varepsilon$, $e_1$, (the edges corresponding to) $\varepsilon'$, and $e_2$,
we actually get a $G$-patch that is not a strip, see Figure \ref{k=2}.
\end{ex}

\begin{co}
\label{2_trees}
Assume there exist disjoint subtrees $G_1$, $G_2$ of the graph $G$ such that 
$G \setminus (G_1 \cup G_2)$ is the disjoint union of $k$ edges, each of which connects $G_1$ to $G_2$.
If $k$ is odd then $G$ is orientably realizable.
\end{co}

\noindent{\sl Proof:} 
The result follows from Theorem \ref{or_cy} and the existence of a $G$-strip obtained by joining the $G_1$-strip to the $G_2$-strip by $k$ switched edge-strips.
\hfill $\Box$

\begin{co}
\label{Petersen}
The Petersen graph is orientably realizable.
\end{co}

\noindent{\sl Proof:}
Deleting the three edges starting from a vertex of the Petersen graph,
we obtain two components: a point and (a graph isomorphic to) the graph in Figure \ref{4cy} (g).
The orientable realizability for the point is trivial, while for the second component two orientable realisations are provided by Theorem \ref{4} (see Figure 12 (m) and (o)). The conclusion follows now by Corollary \ref{conn}.
\hfill $\Box$

\begin{co}
\label{N-or_crit-2}
Assume the graph $G$ can be represented as the union of two disjoint subgraphs, joined by a bridge.
If one of the subgraphs is not orientably realizable
--in particular if it is the graph in Figure 5-- then neither is $G$.
\end{co}

\noindent {\bf Acknowledgement } C. V\^\i lcu was partially supported by the grant PN II Idei 1187 of the Romanian Government.




\bigskip

Jin-ichi Itoh

\noindent {\small Faculty of Education, Kumamoto
University
\\Kumamoto 860-8555, JAPAN
\\j-itoh@gpo.kumamoto-u.ac.jp}

\medskip

Costin V\^\i lcu

\noindent {\small Institute of Mathematics ``Simion Stoilow'' of the
Romanian Academy
\\P.O. Box 1-764, Bucharest 014700, ROMANIA
\\Costin.Vilcu@imar.ro}

\end{document}